\documentclass{article}
\usepackage{amsmath}
\usepackage{amsfonts}
\usepackage{amssymb}
\usepackage{tabularx}
\usepackage{color}
\usepackage{fancyhdr}
\usepackage[all]{xy}
\usepackage{enumerate}

\def\theoremparent{section}
\usepackage{notes}

\begin{document}

\title{Shifted Poisson and Batalin-Vilkovisky structures on the derived variety of complexes}
\author{Slava Pimenov}
\date{\today}
\titlepage
\maketitle


\def\B{\mathcal{B}}
\def\b{\mathfrak{b}}
\def\g{\mathfrak{g}}
\def\h{\mathfrak{h}}
\def\k{\mathfrak{k}}
\def\l{\mathfrak{l}}
\def\n{\mathfrak{n}}
\def\p{\mathfrak{p}}
\def\q{\mathfrak{q}}
\def\ls#1{(\!(#1)\!)}

\def\Ad{\mathop{Ad}\nolimits}
\def\Com{\mathrm{Com}}
\def\Der{\mathop{\mathrm{Der}}\nolimits}
\def\diag{\mathrm{diag}}
\def\Map{\mathrm{Map}}
\def\PCom{\mathrm{PCom}}
\def\per{\mathrm{per}}
\def\RCom{\mathrm{RCom}}
\def\RPCom{\mathrm{RPCom}}
\def\tr{\mathrm{tr}}

\def\sectionseparator{\vskip 5em}
\setcounter{section}{-1}
\section{Introduction.}

Let $V^\bullet$ be a graded vector space of finite total dimension over field $k$. The derived variety of complexes (\cite{KP}) is a dg-scheme
$\RCom(V^\bullet)$ characterized by the following property: for a commutative differential graded algebra $A$ over $k$ the set of maps
$\Hom(\Spec A, \RCom(V^\bullet))$ is the set of differentials in the total complex of $A \tensor V^\bullet$ ``going in the direction of $V^\bullet$'',
making it a twisted complex.

In this paper we will describe a 1-shifted Poisson structure on the quotient of $\RCom(V^\bullet)$ by the infinitesimal action of a subgroup of
automorphisms of $V^\bullet$ with superdeterminant $1$. 
Moreover, this Poisson structure can be upgraded (up to a quasi-isomorphism) to a structure of Batalin-Vilkovisky manifold. Since the ring of functions
on $\RCom(V^\bullet)$ can be explicitly expressed as a cochain complex of certain graded Lie algebra associated to $V^\bullet$ we look at this
Lie algebra for a source of such structures.

\begin{nparagraph}
If $G$ is a classical Poisson-Lie group, the Schouten-Nijenhuis bracket on the polyvector fields translates to a bracket on the differential
forms $\Omega^\bullet G$. Together with the de Rham differential and product it equips $\Omega^\bullet$ with a differential Gerstenhaber algebra
structure (\cite{Koszul}).

We mimic this construction in the case of $n$-shifted graded Lie bialgebra $\g$ to obtain an $(n+1)$-shifted Poisson structure on the cochain complex
$C^\bullet(\g, k)$.
\end{nparagraph}

\begin{nparagraph}
In the second section we apply this result in the case of derived variety of complexes and using the standard bialgebra sturcture on
$\gl(m, n)$ obtain a 1-Poisson structure on the infinitesimal quotient
$\RCom(V^\bullet) / \l_1$, where $\l_1$ is a subalgebra of $\gl(V^\bullet)$ of supertrace $0$ elements.

We also define a derived variety of 1-periodic complexes, that possesses a non-shifted Poisson structure,
extending Kirillov-Kostant structure on $\gl_n^*$.
\end{nparagraph}

\begin{nparagraph}
Although the cochain complex of a graded Lie bialgebra does not necessarily possess a Batalin-Vilkovisky operator, in some cases it is quasi-isomorphic
as a differential Gerstenhaber algebra to a differential BV-algebra. In the last section we give a sufficient condition for existence of such
quasi-isomorphism, which applies to the case of the quotient $\RCom(V^\bullet) / \l_1$.
\end{nparagraph}

The author would like to thank Kavli IPMU for providing wonderful working conditions, and also Mikhail Kapranov and Alexander Voronov for their
invaluable advice in preparation of this paper.

\sectionseparator
\section{Shifted Lie bialgebras.}

Let $V^\bullet = \bigoplus V^i$ be a graded vector space over field $k$, such that each $V^i$ is of finite dimension.

\begin{definition}
An $n$-shifted Lie bialgebra is a graded Lie algebra $(V^\bullet, [-,-])$, with a cobracket
$\phi\from V^\bullet[-n] \to V^\bullet[-n] \tensor V^\bullet[-n]$,
such that the dual map $\phi^*$ is a graded Lie algebra structure on $(V^\bullet[-n])^*$, and the
composition
$$
\xymatrix{
V \ar^-{\phi[n]}[r] & (V[-n] \tensor V[-n])[n] \ar^-\isom[r] & (V \tensor V)[-n]
}
$$
is a $1$-cocycle with coefficients in $(V \tensor V)[-n]$, i.e., it satisfies the following condition:
$$
\phi([x, y]) = \ad_x\phi(y) - (-1)^{\bar{x}\bar{y}} \ad_y\phi(x),
$$
where $\bar{x}$ and $\bar{y}$ are degrees of $x$ and $y$ in $V^\bullet$.
\end{definition}

Abusing notation we denote the cocycle also by $\phi$.

\begin{definition}
An $n$-shifted Manin triple is a triple of graded Lie algebras $(\p, \p_+, \p_-)$, with non-degenerate invariant symmetric bilinear form
$(-,-)\from S^2 \p \to k[n]$, such that $\p = \p_+ \oplus \p_-$ as a graded vector space, $\p_+$ and $\p_-$ are Lie subalgebras of $\p$ isotropic
with respect to the bilinear form.
\end{definition}

We have the following analog of the classical result:
\begin{lemma}
$n$-shifted Lie bialgebra structures on $\g$ are in bijection with $n$-shifted Manin triples $(\p, \p_+, \p_-)$ with $\p_+ = \g$.
\end{lemma}

We begin with a Manin triple and will construct the bialgebra structure on $\g = \p_+$. Using the bilinear form, we can identify $\p_-$ with $\p_+^*[n]$ by sending $x$ to $(x, -)$. Restriction of the Lie bracket to $\p_-$ defines the cobracket $\phi^*$.
Let $\{e_i\}$ be a basis of $\g$, and $\{f^i\}$~--- the dual basis of $\g^*$. We will denote the element of $\p_-$ corresponding to
$f^i$ by $f^i[n]$, so that $(f^i[n], e_j) = \delta^i_j$.

Let $[e_i, e_j] = c_{ij}^k e_k$ and $\phi^*(f^i[n], f^j[n]) = \gamma^{ij}_k f^k[n]$. Using the invariance of the bilinear form we find that
$$
[e_j, f^i[n]] = (-1)^{\wbar{e_i}(n+1)} \gamma^{ik}_j e_k - (-1)^{\wbar{e_j}(\wbar{f^i} - n)} c_{jk}^i f^k[n].
$$

Expanding the Jacobi identity we obtain
\begin{multline*}
[[e_i, e_j], f^k[n]] - [e_i, [e_j, f^k[n]]] + (-1)^{\wbar{e_i}\wbar{e_j}}[e_j, [e_i, f^k[n]]] = \\
(-1)^{1+(\wbar{f^k} - n)\wbar{e_m}} c_{mp}^k c_{ij}^m f^p[n] + (-1)^{\wbar{e_k}(n + 1)} \gamma^{kq}_m c_{ij}^m e_q + \\
(-1)^{1+\wbar{e_k}(n+1)} \gamma^{km}_j c_{im}^q e_q + (-1)^{\wbar{e_j}(\wbar{f^k} - n) + \wbar{e_m}(n+1)} \gamma^{mq}_i c_{jm}^k e_q + \\
(-1)^{1+\wbar{e_j}(\wbar{f^k} - n) + \wbar{e_i}(\wbar{f^m} - n)} c_{jm}^k c_{ip}^m f^p[n] + \\
(-1)^{\wbar{e_i}\wbar{e_j} + \wbar{e_k}(n+1)} \gamma^{km}_i c_{jm}^q e_q +
(-1)^{1+\wbar{e_i}\wbar{e_j} + \wbar{e_i}(\wbar{f^k}-n) + \wbar{e_m}(n+1)} \gamma^{mq}_j c_{im}^k e_q + \\
(-1)^{1+\wbar{e_i}\wbar{e_j} + \wbar{e_i}(\wbar{f^k}-n) + \wbar{e_j}(\wbar{f^m}-n)} c_{im}^k c_{jp}^m f^p[n]
\end{multline*}

Combining coefficients of $f^p[n]$ we obtain the Jacobi identity for $\g$. Now we want to show that the coefficients of $e_q$ provide the
cocycle condition for $\phi$. Using structure constants $\gamma$ we write
$$
\phi(e_k) = (-1)^{\wbar{e_i}\wbar{e_j} + n\wbar{e_i}} \gamma^{ij}_k e_i \tensor e_j [-n] \in (\g \tensor \g)[-n],
$$
and the cocycle condition becomes
\begin{multline*}
(-1)^{\wbar{e_k}\wbar{e_q} + n\wbar{e_k}} \gamma^{kq}_m c_{ij}^m = (-1)^{\wbar{e_m}\wbar{e_q} + n\wbar{e_m}} \gamma^{mq}_j c_{im}^k +
(-1)^{\wbar{e_k}\wbar{e_m} + n\wbar{e_k} + \wbar{e_i}\wbar{e_k}} \gamma^{km}_j c_{im}^q +\\
(-1)^{1+\wbar{e_i}\wbar{e_j}+\wbar{e_m}\wbar{e_q} + n\wbar{e_m}} \gamma^{mq}_i c_{jm}^k +
(-1)^{1+\wbar{e_i}\wbar{e_j}+\wbar{e_k}\wbar{e_m} + n\wbar{e_k} + \wbar{e_j}\wbar{e_k}} \gamma^{km}_i c_{jm}^q.
\end{multline*}
It remains to compare the signs of the corresponding coefficients.

The proof of the lemma in the other direction is similar.

\begin{proposition}
For an $n$-shifted Lie bialgebra $\g$ the dual $\g' = \g^*[n]$ is also an $n$-shifted Lie bialgebra.
\end{proposition}
Immediately follows by interchanging roles of $\p_+$ and $\p_-$ in the Manin triple associated to $\g$.

\begin{proposition}
\label{diff-poisson}
Let $(\g, [-,-], \phi)$ be an $n$-shifted Lie bialgebra, then the cochain complex with trivial coefficients
$C^\bullet(\g) = C^\bullet(\g, k) = S^\bullet(\g^*[-1])$ is a differential
$(n+1)$-shifted Poisson algebra.
\end{proposition}

We will denote by $|a|$ the degree of a cochain $a$ in $C^\bullet(\g)$, so that for any $x \in \g^*$, considered as a cochain we have
$|x| = \bar{x} + 1$. Define the Poisson bracket on $\g^*[-1]$ as
$\<x[-1], y[-1]\> = \phi^*(x[n], y[n]) \in (\g^*)^{\bar{x} + \bar{y} - n} \subset C^{|x| + |y| - n - 1}(\g)$, and extend to
the entire cochain complex using the derivation relation $\<x, yz\> = \<x, y\>z + (-1)^{(|x| - n - 1)|y|} y\<x, z\>$.
The commutation and Jacobi identities for the Poisson bracket are equivalent to those of the cobracket $\phi^*$.

It remains to show the compatibility of the derivation $d$ of the cochain complex and the Poisson bracket:
\begin{equation}
\label{diff-poisson-bracket}
d\<x, y\> = \<dx, y\> + (-1)^{|x| - n - 1} \<x, dy\>.
\end{equation}
Using notations from the previous lemma we have
\begin{align*}
\<f^i[-1], f^i[-1]\> &= \gamma^{ij}_k f^k[-1], \\
d(f^k[-1]) &= (-1)^{\wbar{e_p}\wbar{e_q} + \wbar{e_p}} c_{pq}^k f^p[-1] f^q[-1],
\end{align*}
and it is easy to see that the relation (\ref{diff-poisson-bracket}) is equivalent to the cocycle condition for the Lie bialgebra.

\begin{corollary}
\label{poisson-cohomology}
Let $(\g, [-,-], \phi)$ be an $n$-shifted Lie bialgebra, then the cohomology ring $H^\bullet(\g, k)$ is an $(n+1)$-shifted Poisson algebra.
\end{corollary}

\begin{example}
Let $\g$ be a graded Lie algebra, it has always the trivial $n$-shifted Lie bialgebra structure given by the zero cobracket. The dual
Lie bialgebra $\g'$ is an abelian Lie algebra with an $n$-shifted cobracket. It's cohomology
$H^\bullet(\g', k) = C^\bullet(\g') = S^\bullet(\g[-n-1])$ is an $(n+1)$-shifted Poisson algebra. For $n = 0$ this is the usual Gerstenhaber
algebra structure on the $\Lambda^\bullet\g$ (in fact a Batalin-Vilkovisky algebra).
\end{example}

\begin{example}
\label{gl-m-n}
Now let $\g$ be the graded Lie algebra $\gl(m, n)$, and fix a Borel subalgebra $\b_+$ in $\g$. Denote by $\b_-$ the opposite Borel subalgebra.
We use Manin triples to define a (non-shifted) Lie bialgebra structure on it. $\g$ is spanned
by elementary matrices $e_{ij}$ with $1$ in position $(i, j)$ and $0$ otherwise, and diagonal matrices $a_i$ with $1$ in position
$(i, i)$. Consider the supertrace bilinear form $(a, b) = \tr(ab)$ on $\g$. In this basis it is given by
$(e_{ij}, e_{ji}) = (-1)^{\alpha_i}$, $(a_i, a_i) = (-1)^{\alpha_i}$ and $0$
otherwise. Here $\alpha_i$ is the degree of the $i$-th basis vector in the defining representation of $\g$, i.e. $\alpha_i = 0$ for
$1 \le i \le m$ and $\alpha_i = 1$ for $m < i \le m + n$. This is a non-degenerate symmetric $\g$-invariant bilinear form.

Consider the triple consisting of $\p = \g \oplus \g$ as a graded Lie algebra, $\p_+$~--- the image of the diagonal map from $\g$ to $\p$,
and $\p_- = \{(x, y) \in \b_+ \oplus \b_- \mid p(x) + p(y) = 0\}$, where $p\from \b_\pm \to \h$ are projections of the Borel subalgebras
on the Cartan subalgebra.
Extend the bilinear form from $\g$ to $\p$ by $((x_1, y_1), (x_2, y_2)) = (x_1, x_2) - (y_1, y_2)$. This turns the triple $(\p, \p_+, \p_-)$
into a Manin triple.

If $\q$ is parabolic subalgebra of $\g$, containing $\b_+$, then the orthogonal $\q^\perp \in \g^*$ is identified via the form on $\p$
with a nilpotent ideal in $\b_-$. Therefore the cobracket $\phi$ restricts to $\q$ and obviously satisfies the cocycle condition.

\end{example}

\begin{example}
\label{ex-nonzero-bracket}
Let $\b \subset \gl_n$ is a Borel subalgebra with the restriction of the standard bialgebra structure on $\gl_n$. In this case the cohomology groups
$H^\bullet(\b, k) = \Lambda^\bullet \h^*$, where $\h$ is the Cartan subalgebra of $\gl_n$, and the induced Poisson bracket is trivial.

However this is not the case for the graded Lie algebras. For example, consider a Borel subalgebra $\b \subset \sl(1, 1)$. It has trivial bracket, and
$H^\bullet(\b, k) = k[x] \tensor \Lambda[t]$, where $x$ is the dual of $e_{12} \in \b$, and $t$ dual of the identity matrix $h$. The restriction of the
standard Lie bialgebra structure of $\gl(1, 1)$ gives the cobracket $\phi(e_{12}) = {1\over 2} h \wedge e_{12}$. The Poisson bracket on $H^\bullet(\b, k)$
is defined on the generators as
$$
\<x, x\> = \<t, t\> = 0,\quad \<t, x\> = x.
$$
\end{example}

\begin{example}
\label{ex-poisson-lie-group}
Consider a connected compact Poisson-Lie group $G$. Denote by $\Omega^\bullet G$ and $P^\bullet G$ the differential forms and polyvector fields
respectively. The Poisson structure is given by a bivector $w \in P^2 G$, such that $[w, w] = 0$. This bivector induces the map
$p\from \Omega^1 G \to P^1 G$ defined by $p(\alpha)(\beta) = w(\alpha, \beta)$. It can be extended to a morphism of complexes
$p\from (\Omega^\bullet G, d) \to (P^\bullet G, d_w = [w, -])$.

Now, let $i_w\from \Omega^n G \to \Omega^{n - 2} G$ denote the contraction with bivector $w$, so that $i_w(\alpha \beta) = w(\alpha, \beta)$.
Define degree $-1$ map $B$ as the commutator $[i_w, d]$. The induced bracket on the forms
$$
\<\omega_1, \omega_2\> = (-1)^{\wbar{\omega_1}}(B(\omega_1 \omega_2) - B(\omega_1) \omega_2 -
(-1)^{\wbar{\omega_1}}\omega_1 B(\omega_2))
$$
is mapped to the Schouten bracket of polyvector fields by $p$. Since $i_w$ provides homotopy for $B$ the induced operator on cohomology and the
corresponding bracket vanish.

Recall that in this case the inclusion $\Lambda^\bullet \g^* \into \Omega^\bullet G$ as left invariant forms is a quasi-isomorphism.
The bialgebra structure on $\g$ is defined by providing a bracket on $\g^*$, namely $\phi^*(x, y) = (d i_w(xy))_e$, where $x$ and $y$
are considered both as elements of $\g^*$ and left invariant forms on $G$, and the subscript denotes taking value at the identity point $e$ in $G$.

Even though the map $B$ itself does not preserve left invariant forms, the induced bracket does, and it coincides with the Lie bialgebra
cobracket. Indeed, since $B = [i_w, d]$, we have
$$
\<\alpha, \beta\> = d i_w (\alpha \beta) - i_w d(\alpha \beta) + i_w(d\alpha) \beta - \alpha i_w (d\beta).
$$
Expanding the second term using $d(\alpha \beta) = (d\alpha)\beta - \alpha d\beta$, we find
$$
\<\alpha, \beta\> = d i_w (\alpha \beta) - i_{p(\beta)} (d\alpha) + i_{p(\alpha)} (d\beta).
$$
For any Poisson-Lie group $w_e = 0$, therefore $\<x, y\>_e = d i_w (xy)_e = \phi^*(x, y)$. Since the inclusion as left invariant forms
induces isomorphism $H^\bullet(\g, k) \isom H^\bullet(G, k)$, we see that the $1$-shifted Poisson structure on $H^\bullet(\g, k)$ is trivial.

For the general facts on Poisson-Lie groups used here, we refer to \cite{LW}, \cite{Koszul}.
\end{example}

\begin{nparagraph}{\bf Relation to the 2-shifted symplectic structure on $BG$.}
\label{lagrangians}

Recall the definition of the 2-shifted symplectic structure on $BG$ from \cite{PTVV}.
Let $G$ be an affine smooth algebraic group over field $k$. Assume that it's Lie algebra $\g$ possesses a non-degenerate invariant symmetric
bilinear form. The cotangent complex on the quotient stack $BG$ is $\g^*[-1]$, so that the bilinear form considered as an element of
$$
H^0_{\mathrm{Hoch}}(G, \Lambda^2 (\g^*[-1])) \isom (S^2 \g^*)^G,
$$
is a 2-shifted 2-form on $BG$.

Now assume that $G$ is a Drinfeld double of a Poisson Lie group $K$. Then we have a Manin triple of corresponding Lie algebras
$(\g, \g_+ = \k, \g_- = \k^*)$, and $\g$ is equipped with a symmetric bilinear form, such that $\k$ is identified with an isotropic
subalgebra of $\g$. Then the map of the quotient stacks $i\from BK \to BG$ is Lagrangian. It is straightforward to see, using the characterization
of Lagrangian maps in \cite{Ca}. The cotangent complex $\mathbb{L}_{BK}$ is identified with $\k^*[-1] \isom \g_-[-1]$, and $T_{BK} = \k[1]$
is clearly the kernel of the map $i^*T_{BG} \to \mathbb{L}_{BK}[2]$, induced by the bilinear form.

Here we work with the profinite completion $\hat G$ of $G$ at identity. The ring of functions on $\hat G$ is $(U\g)^*$, and for the 
quotient $B\hat G$ the Hochschild cochain complex $C_{\mathrm{Hoch}}^\bullet(U\g, k)$ is quasi-isomorphic to the cochain complex of the
Lie algebra $C^\bullet(\g, k)$. The results of this section can be viewed as construction of a $1$-shifted Poisson structure on $B\hat K$.

It fits well into the point of view of an $n$-Poisson structure as a bivector of degree $-n$ on $B\hat K$ (\cite{PTVV}). Since the tangent bundle to
$\hat K$ is $\k[1]$, for $n = 1$ we have
$$
\H^0(C^\bullet_{\mathrm{Hoch}}(\hat K, S^2 (T_{B\hat K}[-n-1]))[n+2]) = H_{\mathrm{Hoch}}^1(\hat K, \Lambda^2 \k) = H^1(\k, \Lambda^2 \k).
$$
Which is a 1-cocycle with coefficients in $\Lambda^2 \k$. The co-Jacobi condition for the cocycle is equivalent to the Poisson
condition for the bivector.

Because of this analogy we will use the following terminology:
\end{nparagraph}

\begin{definition}
\label{def-inf-quotient}
The infinitesimal quotient of an affine dg-scheme $X = \Spec A$ by an action of a Lie algebra $\g$ is the spectrum of dg-algebra
(in general with components of both positive and negative degrees) $C^\bullet(\g, A)$.
\end{definition}

\begin{nparagraph}
{\bf Formal classifying space of a graded Lie algebra.}
\label{formal-Bg}

In the case when $\g$ is a graded Lie algebra we define the {\em formal classifying space} $B\g$ of $\g$ as affine dg-scheme
$\Spec C^\bullet(\g, k)$. The tangent bundle to $B\g$ is $\g[1]$ with adjoint action of $\g$,
so that the complex of derivations $\Der(C^\bullet(\g))$ is
isomorphic to $C^\bullet(\g, \g[1])$.

Now, if $\g$ is equipped with an $n$-shifted symmetric invariant non-degenerate scalar product $(-,-)$, then $B\g$ is an $(n+2)$-shifted
symplectic manifold with the form
$$
\omega \in H^{n+2}(\g, \Lambda^2 (\g^*[-1])) = H^0(\g, (S^2 \g^*)[n]),
$$
corresponding to the scalar product. Let us describe the Poisson bracket for this symplectic structure.

As before, for an element $a$ of $\g$ we denote by $\wbar a$ it's degree in $\g$. Let $p\from \g \to \g^*[n]$ be the isomorphism
induced by the scalar product, then the degree of $p(a)$ as an element in $C^\bullet(\g)$ is $|p(a)| = \wbar a + n + 1$.
On the linear functions the $(n+2)$-shifted Poisson bracket is given by
$$
\{p(a), p(b)\} = (-1)^{\wbar a} (a, b).
$$
The Jacobi identity is trivially satisfied. To establish the compatibility with the differential:
$$
\{dp(a), p(b)\} + (-1)^{|p(a)| - n} \{p(a), dp(b)\} = d\{p(a), p(b)\} = 0,
$$
we observe that $\{p(a), dp(b)\} = \Ad_a (p(b))$, then use $\g$-invariance of isomorphism $p$ and
commutation of the Lie bracket of $\g$.

Now let $(\g, \g_+ = \k, \g_- = \k^*[n])$ be an $n$-shifted Manin triple. Denote by $i\from B\k \to B\g$ the map of the formal classifying spaces induced
by the inclusion of $\k$ as a subalgebra of $\g$. We have the short exact sequence of $\k$-modules $\k[1] \to \g[1] \to \g_-[1] \isom \k^*[n+1]$.
It identifies the tangent bundle $T_{B\k} = \k[1]$ with the kernel of the map $i^* T_{B\g} \to L_{B\k}[n+2]$. Therefore the map
$i\from B\k \to B\g$ is Lagrangian (\cite{Ca}).

\end{nparagraph}

\sectionseparator
\section{Poisson structures on the derived varieties of complexes.}
\label{sec-der-complex}

Let $V^\bullet$ be a graded vector space over a field $k$ of characteristic $0$, with only finitely many non-zero spaces, and each $V^i$ of finite dimension.
Recall construction of $\RCom(V^\bullet)$ from \cite{KP}.
Consider the graded Lie algebra $\End(V)$ of endomorphisms of the total vector space or all degrees. Denote by $\q$ the parabolic subalgebra
$$
\q = \bigoplus_{i \le j} \Hom(V^i, V^j)[i - j].
$$

We have the Levi decomposition $\q = \l \oplus \n$, where $\l$ is the Levi subalgebra
and $\n$ the nilpotent radical of $\q$:
$$
\l = \gl(V^\bullet) = \bigoplus_i \gl(V^i),\quad \n = \bigoplus_{i > j} \Hom(V^i, V^j)[i - j].
$$

\begin{definition}
\label{der-complex}
The derived variety of complexes $\RCom(V^\bullet)$ is an affine dg-scheme $\Spec C^\bullet(\n, k)$.
\end{definition}

Choose a basis $\{v_i\}$ of $V^\bullet$, and set $\alpha_i = j$ if the corresponding basis vector $v_i \in V^j$. For convenience enumerate
basis vectors so that $\alpha_i \le \alpha_{i'}$ if $i < i'$. Denote by $e_{ij} = v_i \tensor v_j^*$, the basis vector of $\End(V)$
of degree $\wbar{e_{ij}} = \alpha_j - \alpha_i$.

\begin{nparagraph}
{\bf Poisson bracket induced by the standard Lie bialgebra structure.}

We may identify $\q$ with the Lie subalgebra of block upper triangular matrices (including
the diagonal blocks) and $\n$ with the block strict upper triangular matrices.
Let $\b_+$ be the Borel subalgebra of upper triangular matrices, then using the same construction as in example \ref{gl-m-n}, we obtain a Lie bialgebra structure on $\q$.

Now let $\l_1 = \l \cap \sl(V)$ and $\q_1 = \q \cap \sl(V)$, where $\sl(V)$ is the subalgebra of $\End(V)$ consisting of elements with
supertrace $0$. The orthogonal of $\q_1$, identified with $\n_- \oplus k$, is an ideal in $\End(V)^*$ with respect to the Lie cobracket, hence the
bialgebra structure restricts to $\q_1$ as well.

By proposition \ref{diff-poisson} we have a $1$-shifted differential Poisson structures on the $C^\bullet(\q, k)$ and $C^\bullet(\q_1, k)$.
First we will show that Poisson bracket on $C^\bullet(\q)$ induces trivial bracket on the cohomology. Consider the Hochschild-Serre
spectral sequence associated
to the quotient $\q / \n$:
$$
E_2^{pq} = H^p(\l, H^q(\n, k)) \Rightarrow H^{p+q}(\q, k).
$$
Since $\l$ is a reductive Lie algebra, it is sufficient to look at submodule of $\l$-invariants in $H^\bullet(\n, k)$.
Because of weight considerations ($\l$ contains the Cartan subalgebra) $H^\bullet(\n)^\l = H^0(\n) = k$, so that the spectral sequence
degenerates and $H^\bullet(\q, k) = H^\bullet(\l, k)$.

After changing base from $k$ to $\CC$, it remains to show that the Poisson bracket on $H^\bullet(\gl_n(\CC))$ is trivial. This can be done
using the argument similar to example \ref{ex-poisson-lie-group}.

However the Poisson bracket induced on $H^\bullet(\q_1, k)$ in general is not trivial, as was seen in example \ref{ex-nonzero-bracket}.
\end{nparagraph}

\begin{theorem}
The infinitesimal quotient (in the sense of \ref{def-inf-quotient}) of the derived variety of complexes $\RCom(V^\bullet) / \l_1$ has structure
of a 1-shifted Poisson manifold.
\end{theorem}

Follows immediately from the discussion above and observation that
$$
\RCom(V^\bullet)/\l_1 := \Spec C^\bullet(\l_1, C^\bullet(\n, k)) \isom \Spec C^\bullet(\q_1, k).
$$

\begin{remark}
Instead of taking the infinitesimal quotient we could consider the quotient stack $\RCom(V^\bullet) / GL(V^\bullet)_1$,
where $GL(V^\bullet) = \prod_i GL(V^i)$, and $GL(V^\bullet)_1$ is the subgroup consisting of elements
$(g_1, \ldots , g_n)$ with
$$
\prod_i (\det g_i)^{(-1)^i} = 1.
$$
One may ask if the 1-Poisson structure of the theorem lifts to a Poisson bivector
on this quotient stack.
\end{remark}

\begin{nparagraph}{\bf $\gl(n, n)$ as a Drinfeld double and Poisson structures on the derived variety of complexes.}
\label{drinfeld-double}

Let $\p = \gl(n, n)$, fix a Borel subalgebra $\b_+$ in $\p$ and define bilinear form $(-, -)$ as in \ref{gl-m-n}.
Fix also a permutation $\theta$ of the set $\{1, 2, \ldots , n\}$, and denote by $\h^\theta_+$ and $\h^\theta_-$
the subspaces of Cartan
subalgebra $\h$ of $\p$ spanned by $\{a_i + a_{n+\theta(i)}\}$, and $\{a_i - a_{n+\theta(i)}\}$ respectively, $1 \le i \le n$.
Let $\g^\theta = \p^\theta_+ = \n_+ \oplus \h^\theta_+$ and $\p^\theta_- = \n_- \oplus \h^\theta_-$.

It is immediate to see that $\p^\theta_+$ and $\p^\theta_-$ are isotropic subalgebras of $\p$, and the bilinear form identifies
$\p^\theta_-$ with the dual of $\p^\theta_+$. Therefore we have a (non-shifted) Manin triple $(\p, \p^\theta_+, \p^\theta_-)$,
providing $\g^\theta$ with a Lie bialgebra structure.

Notice that this is not a restriction of the standard Lie bialgebra structure on $\gl(n, n)$ to $\g^\theta$, except when $n = 1$,
because in general the standard Lie cobracket doesn't restrict to $\g^\theta$.
\end{nparagraph}

\begin{nparagraph}
We may relate the cohomology with trivial coefficients of $\g^\theta$ to the cohomology of the nilpotent radical $\n_+$ computed in \cite{KP}.
First, let us restrict to the quasi-isomorphic subcomplex $C^\bullet(\g^\theta)^{\h^\theta_+}$ of $\h^\theta_+$-invariant elements.
Now, since
$$
C^\bullet(\g^\theta)^{\h^\theta_+} = \Lambda^\bullet(\h^\theta_+)^* \tensor C^\bullet(\n_+)^{\h^\theta_+},
$$
with the differential induced by that of the cochain complex of $\n_+$ and zero differential in the exterior powers, we find
$$
H^\bullet(\g^\theta) = \Lambda^\bullet(\h^\theta_+)^* \tensor H^\bullet(\n_+)^{\h^\theta_+}.
$$
The $\h^\theta_+$-invariant subspace of $H^\bullet(\n_+)$ is the subspace with weights belonging to the configuration
sector $\mathfrak S(\theta, \b_+)$ (see {\it op. cit.} 2.3).
\end{nparagraph}

\begin{nparagraph}
\label{borel-complex}
Let $V^\bullet$ be a graded vector space with each $V^i$ of dimension at most $1$, such that both total odd and total even dimensions are $n$.
Then after appropriate shuffle permutation $\End(V)$ as $\Z / 2$-graded Lie algebra is isomorphic to $\gl(n, n)$, and ideal $\n$ defined at
the beginning of this section corresponds to the nilpotent radical of some Borel subalgebra $\b_+$. In this setting we have
\end{nparagraph}

\begin{proposition}
For each permutation $\theta \in S_n$, the infinitesimal quotient of $\RCom(V^\bullet)$ by the action of $\h^\theta_+$ has 1-shifted Poisson structure.
Each such quotient is a Lagrangian in the formal classifying space (see \ref{formal-Bg}):
$$
B\gl(V) := \Spec C^\bullet(\gl(V), k).
$$
\end{proposition}
Moreover, in section \ref{sec-BV} we will see that these quotients are up to homotopy BV-manifolds.

\begin{nparagraph}
{\bf Derived variety of 1-periodic complexes.}

As another application of the machinery developed here we would like to give a related example of non-shifted Poisson structure
on the derived variety of 1-periodic complexes.

Let $V^\bullet$ be a graded vector space with $W$ in each degree $i \in \Z$. Define the space $\PCom(W)$ of $1$-periodic complexes on $V^\bullet$
as a subspace of $\End(W)$ consisting of all square-zero elements $d$, so that
$$
\xymatrix{
\ldots \ar^d[r] & W \ar^d[r] & W \ar^d[r] & W \ar^d[r]&  \ldots
}
$$
is an unbounded complex. Let $\{e_{ij}\}$ be a basis of $\End(W)$, and $\{x_{ij}\}$ the corresponding dual basis of $\End(W)^*$. The
Lie bracket on this dual space provides the Kirillov-Kostant Poisson bracket on $\End(W)$.

$\PCom(W)$ is an affine variety, it's ring of functions is the quotient $S^\bullet (\End(W)^*) / I$,
where $I$ is the ideal generated by $\sum_k x_{ik} x_{kj}$, for all $i, j$. It is straightforward
to check that the Kirillov-Kostant bracket induces a well-defined Poisson structure on $\PCom(W)$. The corresponding symplectic foliation
is the decomposition of $\PCom(W)$ into the union of coadjoint orbits.

Before we define the derived version of $\PCom(W)$ consider the next example.
\end{nparagraph}

\begin{example}
\label{series-lie-alg}
Let $A$ be a finite dimensional symmetric Frobenius algebra (in degree $0$), i.e., an associative algebra equipped with an isomorphism of $A$-bimodules $\theta\from A \to \Hom(A, k)$. This is equivalent to having a non-degenerate
symmetric bilinear form on $A$, satisfying the invariance condition $(ab, c)_A = (a, bc)_A$.

Consider the graded Lie algebra $\wtilde A = t A[t]$, with $\deg t = 1$, so that the Lie bracket is defined as $[at^m, bt^n] = (ab - (-1)^{mn} ba) t^{m+n}$.
We want to show that $\wtilde A$ has a $(-1)$-shifted Lie bialgebra structure. In order to do this consider the triple
$(A\ls{t^{-1}}, \wtilde A, A[[t^{-1}]])$, with the bilinear form $(f, g) = \res_\infty(f, g)_A$, i.e., the coefficient of $t^1$ in the series
$(f, g)_A$. This is clearly a non-degenerate bilinear form on $A\ls{t^{-1}}$, it is symmetric since $(-,-)_A$ is symmetric and the form is
non-trivial only if one of the homogenous terms is of even degree.

It remains to show invariance. Let $m, n, l$ are such that $m + n + l = 1$, then using symmetry and invariance of $(-,-)_A$ and the fact that
$(-1)^{mn} = (-1)^{(n + l + 1)n} = (-1)^{nl}$ we have
\begin{multline*}
([at^m, bt^n], ct^l) = (ab - (-1)^{mn} ba, c)_A t =\\
((a, bc)_A - (-1)^{mn}(ba, c)_A)t = (a, bc - (-1)^{nl} cb)_A t = (at^m, [bt^n, ct^l]).
\end{multline*}

The Lie subalgebras $\wtilde A$ and $A[[t^{-1}]]$ are isotropic with respect to this form, so we have a $(-1)$-shifted Manin triple.
By proposition \ref{diff-poisson} the cochain complex $C^\bullet(\wtilde A, k)$ is a differential Poisson algebra.
\end{example}

Now return to the space of 1-periodic complexes. Let $A = \End(W)$, with the trace scalar product $(a, b) = \tr(ab)$, for any
$a, b \in A$. Construct Lie
algebras $\wtilde A$ and $A[[t^{-1}]]$ as in example \ref{series-lie-alg}.

\begin{definition}
The derived variety of 1-periodic complexes on a vector space $W$ is an infinite dimensional affine dg-scheme
$$
\RPCom(W) = \Spec C^\bullet(\wtilde A, k) = \Spec S^\bullet(A[[t^{-1}]]).
$$
\end{definition}
In particular it is easy to see that the classical variety $\PCom(W) = \Spec H^0(\wtilde A, k)$.

From the discussion of example \ref{series-lie-alg} it immediately follows that

\begin{theorem}
The derived variety of 1-periodic complexes $\RPCom(W)$ is a Poisson dg-manifold, and the restriction of the Poisson structure
to $H^0(\wtilde A, k)$ coincides with the standard Kirillov-Kostant Poisson structure described above.
\end{theorem}

\sectionseparator
\section{Batalin-Vilkovisky algebra structure.}
\label{sec-BV}

A BV${}_\infty$-algebra $A$ (cf. \cite{CL}) is a commutative graded algebra equipped with an operator $D\from A \to A[[h]]$
of degree $1$, such that $D^2 = 0$, $D(1) = 0$ and in the expansion
$$
D = \sum_{i = 1}^\infty D_i h^{i-1},
$$
each $D_i$ is a differential operator of order $\le i$.

A differential BV algebra is a BV${}_\infty$-algebra with all $D_i = 0$ for $i \ge 3$. In this case we take $h$ to be of degree 2,
and write $D = d + Bh$, so that $d$ is a differential in $A$ of degree $+1$ and $B$ is an operator of order $\le 2$ and degree $-1$.
It will be convenient for us to unpack this definition.

\begin{definition}
A differential Batalin-Vilkovisky algebra is a mixed complex $(C^\bullet, d, B)$ (in the sense of \cite{L}), such that $(C^\bullet, d)$
is a differential graded algebra and $(C^\bullet, B)$ is a BV-algebra.
\end{definition}

\begin{nparagraph}
Let $\g$ be a non-shifted graded Lie bialgebra, then as we saw before the cochain complex $C^\bullet(\g, k) = S^\bullet(\g^*[-1])$
has a structure of Gerstenhaber algebra. As a graded vector space it can be identified with the chain complex of the dual Lie
algebra $C_\bullet(\g^*, k)$. The chain differential $B$ considered as a degree $-1$ operator on the cochain complex
is a second order differential operator, with $B^2 = 0$ and $B(1) = 0$. It equips $C^\bullet(\g)$ with a BV algebra structure,
generating the Gerstenhaber bracket constructed before: for any $x, y \in C^\bullet(\g)$ of degrees $|x|$ and $|y|$
$$
\<x, y\> = (-1)^{|x|}(B(xy) - (Bx) y - (-1)^{|x|}x By).
$$

Notice that in general it is not a differential BV-algebra, since operators $B$ and $d$ do not commute. However, we will describe a class
of Lie bialgebras such that $C^\bullet(\g, k)$ is a differential BV-algebra at least up to a quasi-isomorphism.

Since $B$ is a second order operator and $d$ is a first order operator, a priori the commutator $[B, d]$ is a second order operator. However
the cocycle condition for the Lie bialgebra ensures that it is in fact a derivation.
\end{nparagraph}

\begin{lemma}
\label{delta-leibniz}
Let $(C^\bullet, d, \<,\>)$ be a differential Gerstenhaber algebra, with the bracket generated by a second order operator $B$.
Then
the commutator $\Delta = Bd + dB \from C^\bullet \to C^\bullet$ is a derivation, i.e. it satisfies the Leibniz identity: $\Delta(xy) = \Delta(x) y + x \Delta(y)$.
\end{lemma}

Since $C^\bullet$ is a differential Gerstenhaber algebra, we have $d\<x, y\> = \<dx, y\> + (-1)^{|x| + 1}\<x, dy\>$. Rewriting brackets using
codifferential $B$, on the left hand side we have
$$
(-1)^{|x|}dB(xy) + (-1)^{|x| + 1} dBx\cdot y + Bx\cdot dy - dx \cdot By + (-1)^{|x| + 1} x\cdot dBy,
$$
and on the right hand side:
\begin{multline*}
(-1)^{|x| + 1}B(dx \cdot y) + (-1)^{|x|} Bdx \cdot y - dx \cdot By - \\
B(x \cdot dy) + Bx \cdot dy + (-1)^{|x|} x \cdot Bdy.
\end{multline*}
Comparing these two expressions we establish the required identity.

\begin{corollary}
Let $(C^\bullet, d, \<,\>, B)$ be as in the lemma, $\Delta = [B, d]$. The kernel $\Ker\Delta$ is a differential BV-algebra, and the inclusion
$\Ker\Delta \into C^\bullet$ is a map of differential Gerstenhaber algebras.
\end{corollary}

Commutation of $\Delta$ with $B$ and $d$ implies that they restrict to the $\Ker\Delta$. Since $\Delta$ is a derivation, $\Ker\Delta$ is also a
subalgebra. Finally the bracket, being generated by $B$, can also be restricted to the kernel.

\begin{corollary}
Let $(C^\bullet, d, B)$ be a differential BV-algebra,
with bracket $\<-,-\>$ induced by $B$. Then $(C^\bullet, d, \<-,-\>)$ is a differential Gerstenhaber algebra.
\end{corollary}

In a mixed complex the commutator $\Delta = 0$, so that the relation in the previous lemma is trivially satisfied. This in turn equivalent to
$d$ being a derivation of the Gerstenhaber bracket.

\begin{corollary}
If $\g$ is a graded Lie bialgebra, then the commutator $\Delta = [B, d]$ is a derivation.
\end{corollary}

\begin{corollary}
Let $\g$ be an involutive graded Lie bialgebra, i.e., the composition
$$
\xymatrix@1{\g \ar^\phi[r] & \Lambda^2 \g \ar^{[,]}[r] & \g}
$$
is trivial. Then $C^\bullet(\g)$, with $d$ and $B$ as before, is a differential BV-algebra.
\end{corollary}

The vanishing of the composition of Lie bracket and cobracket implies $\Delta(x) = 0$ for all $x \in C^1(\g)$. Since $C^\bullet(\g)$ is generated
as an algebra by $C^1(\g)$ the statement follows.

\begin{example}
Let $\g$ be a coboundary Lie bialgebra, i.e., with the cobracket given by $\phi(x) = [r, x]$ for some $r \in S^2(\g[1])$, with
$[r, r] \in (S^3(\g[1]))^\g$. Here $[-,-]$ denotes the Schouten bracket, generated by the differential $\delta$ of the chain complex.
We want to find a condition for $r$ for the cochain complex $C^\bullet(\g)$ to be a mixed complex. It will be more convenient to look
at the dual to the commutator:
$$
\Delta^*(x) = [\delta, \phi](x) = \delta[r, x] + [r, \delta x] = [\delta r, x].
$$

Therefore $\g$ is an involutive Lie bialgebra, i.e., $\Delta = 0$ if and only if $\delta r$ lies in the center of $\g$. For example, this is the case for $\sl(1, 1)$ with the
standard Lie bialgebra structure, where $r = {1\over 2} e_{12} \wedge e_{21}$, and $\delta r$ is the identity matrix.
\end{example}

\begin{example}
\label{standard-r}
Now look at the graded Lie algebra $\g = \gl(m, n)$ with the standard Lie bialgebra structure. It is a coboundary Lie bialgebra with
$$
r = {1\over 2}\sum_{i < j} e_{ij} \wedge e_{ji}.
$$
Although in this case $\delta r$ doesn't in general belong to the center of $\g$, it lies in the Cartan subalgebra of $\g$.
In the case of $\gl(n)$ it's the half-sum of the positive roots: $\delta r = {1\over 2} \diag(n-1, n-3, \ldots, -n+1)$,
and for $\gl(m, n)$, we have $\delta r = {1\over 2} \diag(m+n-1, \ldots, -m+n+1, m+n-1, \ldots, m-n+1)$.
Therefore $\Delta^* = [\delta r, -]$ acts semi-simply on the chain complex. This is an important case due to the following lemma.
\end{example}

\begin{lemma}
\label{delta-semisimple}
Let $\g$ be a graded Lie bialgebra. As before, $d$ and $B$ are the differential and BV operator on the cochain complex $C^\bullet(\g)$.
If $\Delta = [B, d]$ acts semi-simply on $\g^*$, then the embedding $\Ker \Delta \into C^\bullet(\g)$ is
a quasi-isomorphism of differential Gerstenhaber algebras.
\end{lemma}

Using the Leibniz identity of lemma \ref{delta-leibniz} we see that $\Delta$ acts semi-simply on all $C^n(\g)$. Since it commutes with
the differential we have the direct sum decomposition of the cochain complex $C^\bullet(\g) = \bigoplus_\lambda C^\bullet(\g)^\lambda$
into eigenspaces of $\Delta$. Now, by definition, $B$ is a homotopy for $\Delta$, so that $\Delta$ acts trivially on the cohomology
groups of $\g$, therefore $\Ker \Delta = C^\bullet(\g)^{\lambda = 0}$ is a quasi-isomorphic subcomplex of $C^\bullet(\g)$.

\begin{example}
Let us return back to the setting of section \ref{sec-der-complex}. The quotient of the derived variety of complexes $\RCom(V)$ on a graded
vector space $V^\bullet$ by the Lie algebra $\l_1$ of infinitesimal automorphisms with supertrace $0$,
is identified with the spectrum of the differential graded algebra $S^\bullet(\q_1^*[-1])$ where $\q_1$ is the
parabolic subalgebra of endomorphisms of the total space $V$  with $0$ supertrace and of non-negative degree. The r-matrix for the standard Lie bialgebra
is given by the same expression as in the example \ref{standard-r}. We express $r$ as an element in $\Lambda^2 \g$ rather than
$S^2(\g[1])$ to simplify signs.

Because $\delta r$ is a diagonal matrix, it acts semisimply on $S^\bullet(\q_1^*[-1])$, and according to the lemma
the kernel of the coadjoint action $\Ad^*_{\delta r}$ is a quasi-isomorphic Gerstenhaber subalgebra. In other words

\begin{theorem}
The infinitesimal quotient of the derived variety of complexes $\RCom(V^\bullet) / \l_1$ is a homotopy BV-manifold.
\end{theorem}
\end{example}

\begin{example}
It is clear from the construction of the Lie bialgebra structure on $\g^\theta$ from section \ref{drinfeld-double} that the action
of the commutator $\Delta$ on $(\g^\theta)^*$ is semi-simple, and according to lemma \ref{delta-semisimple}, $\Ker \Delta$ is a quasi-isomorphic
Gerstenhaber subalgebra of $C^\bullet(\g^\theta)$.

\begin{proposition}
Under assumptions of example \ref{borel-complex} the infinitesimal quotient $\RCom(V^\bullet) / \h^\theta_+$ is a homotopy BV-manifold.
\end{proposition}
\end{example}

\vfill\eject

\end{document}